\documentclass[11pt]{article}
\usepackage{ifthen}
\usepackage{latexsym}
\usepackage{graphicx}
\usepackage{amsmath,mathrsfs}
\usepackage{amssymb}
\usepackage{amsfonts}
\usepackage{wasysym}
\usepackage{amsthm,amssymb,amscd}
\usepackage{epsfig}
\usepackage{epstopdf,fancyhdr}
\usepackage{extarrows,caption,subfigure}
\usepackage{amsmath}
\usepackage{amsfonts}
\usepackage{amssymb}
\usepackage{graphics}
\usepackage{mathrsfs,fancyhdr,authblk}

\theoremstyle{plain}
\newtheorem{thm}{Theorem}[section]
\newtheorem{prop}[thm]{Proposition}
\newtheorem{lem}[thm]{Lemma}
\newtheorem{cor}[thm]{Corollary}

\theoremstyle{definition}

\newtheorem{rem}{Remark}
\newtheorem{defn}{Definition}[section]
\newtheorem{eg}[thm]{Example}

\numberwithin{equation}{section}

\newcommand{\bthm}{\begin{thm}}
\newcommand{\ethm}{\end{thm}}
\newcommand{\bprop}{\begin{prop}}
\newcommand{\eprop}{\end{prop}}
\newcommand{\bcor}{\begin{cor}}
\newcommand{\ecor}{\end{cor}}
\newcommand{\blem}{\begin{lem}}
\newcommand{\elem}{\end{lem}}
\newcommand{\bca}{\begin{cases}}
\newcommand{\eca}{\end{cases}}
\newcommand{\brem}{\begin{rem}}
\newcommand{\erem}{\end{rem}}
\newcommand{\bpm}{\begin{pmatrix}}
\newcommand{\epm}{\end{pmatrix}}
\newcommand{\bdefn}{\begin{defn}}
\newcommand{\edefn}{\end{defn}}
\newcommand{\bsub}{\begin{subtitle}}
\newcommand{\esub}{\end{subtitle}}
\newcommand{\ben}{\begin{enumerate}}
\newcommand{\een}{\end{enumerate}}
\newcommand{\beg}{\begin{eg}}
\newcommand{\eeg}{\end{eg}}
\newcommand{\beq}{\begin{equation}}
\newcommand{\eeq}{\end{equation}}

\def\ms{\medskip}

\def\dim{{\rm dim\/}}
\def\Aut{{\rm Aut\/}}
\def\csch{{\rm csch\/}}

\def \a {\alpha}

\def \l {\lambda}

\def\R{\mathbb{R} }

\DeclareMathOperator{\SL}{\mathrm{SL}}

\begin{document}

\title{\bf The associated families of semi-homogeneous complete hyperbolic affine spheres}

\author[1]{Zhicheng Lin}
\author[2]{Erxiao Wang}
\affil[1]{Wuhan Institute of Physics and Mathematics, CAS, Wuhan, 430071, China~~~~~~flyriverms@qq.com}
\affil[2]{Department of Mathematics, Hong Kong University of Science and
Technology, Clear Water Bay, Kowloon, Hong Kong ~~~~~~maexwang@ust.hk}

\date{}\maketitle

\noindent{\bf Abstract}
Hildebrand classified all semi-homogeneous cones in $\mathbb{R}^3$ and computed their corresponding complete hyperbolic affine spheres. We compute isothermal parametrizations for Hildebrand's new examples. After giving their affine metrics and affine cubic forms,  we construct the whole associated family for each of Hildebrand's examples. The generic member of these affine spheres is given by Weierstrass $\mathscr{P}, ~\zeta ~\text{and} ~\sigma$ functions. In general any regular convex cone in $\R^3$ has a natural associated $S^1$-family of such cones, which deserve further studies.
\bigbreak

\noindent {\bf Keywords} hyperbolic affine spheres, isothermal coordinates, Weierstrass elliptic functions, Monge-Amp\`{e}re equation, Tzitz\'{e}ica equation

\vskip 1mm
\noindent{\bf 2000 MR Subject Classification} 37Kxx; 53A15

\section{Introduction}

Classical equiaffine differential geometry investigates the properties of hypersurface $r(x_1,x_2,\cdots,x_n)$ in $\R^{n+1}$ invariant under the equiaffine transformations $r \to Ar+v$, where $A \in \SL_{n+1}(\R)$ and $ v \in \R^{n+1}$. Although the Euclidean angle is no longer invariant, there exists an affine invariant transversal vector field along $r$, called the affine normal.
The affine metric and the affine cubic form give a set of complete affine invariants. The fundamental theorem tells us that any given affine metric and affine cubic form satisfying certain compatibility conditions (Gauss-Codazzi type) will determine an affine hypersurface uniquely up to equiaffine transformations (see \cite{Lin-Wang-Wang}). The affine metric (conformal to the Euclidean second fundamental form) is definite if and only if the hypersurface is locally strictly convex. The simplest interesting class of hypersurfaces is the  affine spheres, which are defined by the condition that all  affine normal lines meet in a point.
Specially, definite affine spheres (with a mean curvature H and with center at the origin or infinity) can be represented as a graph of a locally strictly convex function $f$ if and only if the Legendre transform $u$ of
$f$ solves a Monge-Amp\`{e}re equation (see Calabi \cite{Eugenio Calabi.{1972}}):
\begin{equation}\label{eqMA}
\det \left(\frac{\partial^2 u}{\partial y_i \partial y_j}\right) = \left\{ \begin{array}{ll}
(H u )^{-n-2}, & \textrm{if $H \neq 0$,}\\
1, & \textrm{if $H=0$.}
\end{array} \right.
\end{equation}
Cheng \& Yau \cite{Cheng-Yau.1986} showed that on a bounded convex domain there is for $H < 0$ a unique negative convex solution of \eqref{eqMA} extending continuously to be $0$ on the boundary. This leads to a beautiful geometric picture of any complete hyperbolic ($H < 0$) affine sphere $r$ (conjectured by Calabi \cite{Eugenio Calabi.{1972}}): it is always asymptotic to the boundary of the cone given by the convex hull of $r$ and its center; and conversely the interior of any regular convex cone is foliated by complete hyperbolic affine spheres asymptotic to it with all $H < 0$ and with centers at the vertex.

However, explicit representations are rarely known except for homogeneous cones.
Hildebrand classified three-dimensional regular convex cone with an automorphism group of dimension at least 2, which he also called \emph{`semi-homogeneous'}. In this case, he reduced the Monge-Amp\`{e}re equation to some ODE and was able to solve it with elliptic integrals in \cite{Roland Hildebrand.{2013}}.

It is known that 2 dimensional affine spheres have additional features than higher dimensional ones. For example, the affine metric can be expressed simply in some natural isothermal coordinates (definite case) or asymptotic coordinates (indefinite case).  In addition, each 2 dimensional affine sphere has a natural associated family of the same type with the same affine metric but different affine cubic forms. In particular, this implies that the structure equations actually form an integrable system, called Tzitz\'{e}ica equation or the affine $\mathfrak{a}_2^{(2)}$-Toda field equation (see \cite{DoEi01}, \cite{Dor-Wang01}, \cite{Lin-Wang-Wang}).  It is then natural to ask for the isothermal parametrizations and the associated families of Hildebrand's new examples. Specially, one wish to see the natural associated family of cones. This paper will answer these questions. 

Our extensive studies of definite affine spheres were motivated by Loftin, Yau \& Zaslow's `trinoid' construction (see \cite{LoYauZa05}) with applications in mirror symmetry.
The dressing actions on proper definite affine spheres and soliton examples have been presented in \cite{Lin-Wang-Wang}.  The Permutability  Theorem and group structure of dressing actions will be presented  in a subsequent  paper \cite{Wang-lin-Wang}. Their Weierstrass or DPW representations have been studied in \cite{Dor-Wang01}, using an Iwasawa decomposition of certain twisted loop group. The equivariant solutions have been constructed in \cite{Dor-Wang02}.

The rest of the paper is organized as follows. In section 2, we introduce the fundamental concepts of affine sphere and review the main results in \cite{Roland Hildebrand.{2013}}. In section 3, we compute the isothermal parametrizations and the affine invariants of Hildebrand's examples case by case. We also get the whole associated family for each case by solving their structure equations. In general any regular convex cone in $\R^3$ has a natural associated $S^1$-family of such cones, which deserve further studies.
\section{Affine spheres and semi-homogeneous cones}

Let
\begin{align*}
& r: M \rightarrow \mathbb{R}^3 \\
&(x,y) \mapsto r(x,y)
\end{align*}
be an immersion with a non-degenerate second fundamental form. Introduce
\[
L = |r_x, r_y, r_{xx}|, ~~ M = |r_x, r_y, r_{xy}|, ~~ N=|r_x, r_y, r_{yy}|,
\]
where $|\cdot,\cdot,\cdot|$ denotes the standard determinant in $ \mathbb{R}^3$. The quadratic form
\[
g = \frac{Ldx^2+2Mdxdy+Ndy^2}{|LN-M^2|^{1/4}}
\]
is equiaffine invariant, and it is called the \textbf{affine} (or \textbf{Blaschke}) \textbf{metric} of the immersion. The surface is said to be \textbf{definite} or \textbf{indefinite} if this metric $g$ is so. The surface is definite if and only if it is locally strictly
convex. A transversal vector field $\xi$ on a surface $r(M)$ is called \textbf{affine normal} if it satisfies
\[
\xi = \frac{1}{2} ~\Delta^g r ,
\]
where $\Delta^g$ is the Laplace-Beltrami operator of the affine metric. The affine metric and the affine normal are then uniquely determined (up to
a sign) by requiring the following decomposition into tangential and transverse
components:
\[
 D_X r_{\ast}Y = r_{\ast}(\nabla_X Y) + g(X,Y)\xi ~,
\]
where $D$ is the canonical flat connection on $\R^3$ and X, Y are any tangent vector fields on the surface. The \textbf{affine cubic form}
measures the difference between the induced Blaschke connection $\nabla$ and $g$'s Levi-Civita connection
$\nabla^g$:
\[
 C(X,Y,Z):= g(\nabla_X Y-\nabla^g_X Y,Z).
\]
It is actually symmetric in all 3 arguments and is a 3rd order invariant.

If all affine normals of $r$ meet in one point (the center), then the surface is called \textbf{affine sphere}. If this point is not infinite it may be chosen as the origin of $\mathbb{R}^3$ so that
\[
\xi =- H r, ~~~~H: M\rightarrow\mathbb{R}.
\]
$H$ is called the \textbf{affine mean curvature}, and such surface is called \textbf{proper affine sphere}. For a convex proper affine sphere, if the center is outside the surface, it is called \textbf{hyperbolic}.

In the sequel we consider the cone $K \subset \mathbb{R}^3$, with the vertex at the origin. The complete hyperbolic affine spheres foliating $K$ are the level sets of the solution $F:K^\circ\rightarrow\mathbb{R}$ of the following Monge-Amp\`{e}re equation (see \cite{Cheng-Yau.1980}, \cite{Cheng-Yau.1982}, \cite{Fox.2012}, \cite{Loftin.2002} for detail):
\begin{equation}\label{eqMA2}
    \det F'' = e^{2F}, ~~~~ F|_{\partial K} = +\infty, ~~~~ F''>0.
\end{equation}

Let $\Aut K$ denote the automorphism group of the cone $K$.  The equiaffine invariance of \eqref{eqMA2} implies F is invariant under unimodular  automorphisms $g \in \Aut K$ (see \cite{Fox.2012} or \cite{Roland Hildebrand.2012}). Then the dimension of \eqref{eqMA2} can be effectively reduced by the generic dimension of the orbits of $\Aut K$. If there are orbits of dimension 2, the PDE \eqref{eqMA2} will reduce to an ODE. This kind of cone is also called \textbf{semi-homogeneous cone}. Hildebrand provided an important classification theorem of  such cones:

\bthm \label{thm}
(\cite{Roland Hildebrand.{2013}},Theorem 3.2) Let $K \subset \mathbb{R}^3$ be a regular convex cone such that \dim \ \Aut $K \geq 2$. Then K is isomorphic to exactly one of the following cones.

1. the cone obtained by the homogenization of the epigraph of the exponential function;

2. the positive orthant $\mathbb{R}^3_{++}$;

3. the cone given by $\{x \, | \,  x_2 \leq x_1^{1/p} x_3^{1/q}, x_1 \geq 0, x_3 \geq 0\}$ for some $p \in [2,\infty)$, $\frac{1}{p}+\frac{1}{q}=1$;

4. the cone given by $\{x \, | \,  -\a x_1^{1/p} x_3^{1/q} \leq x_2 \leq x_1^{1/p} x_3^{1/q}, x_1 \geq 0, x_3 \geq 0\}$ for some $p \in [2,\infty)$, $\frac{1}{p}+\frac{1}{q}=1$, $\a \in (0,1]$;

5. the cone given by $\{x \, | \,  0 \leq x_2 \leq x_1^{1/p} x_3^{1/q}, x_1 \geq 0, x_3 \geq 0\}$ for some $p \in [2,\infty)$, $\frac{1}{p}+\frac{1}{q}=1$.
\ethm

We choose $ p = 5$ for each case and give the pictures of these five cases of semi-homogeneous cones in Figure 1. It could be useful for knowing the shape of semi-homogeneous cones.

\begin{figure}[!htb]
\centering
\subfigure{
\begin{minipage}[t]{0.3\textwidth}
\centering
\includegraphics[width=1.0 \columnwidth]{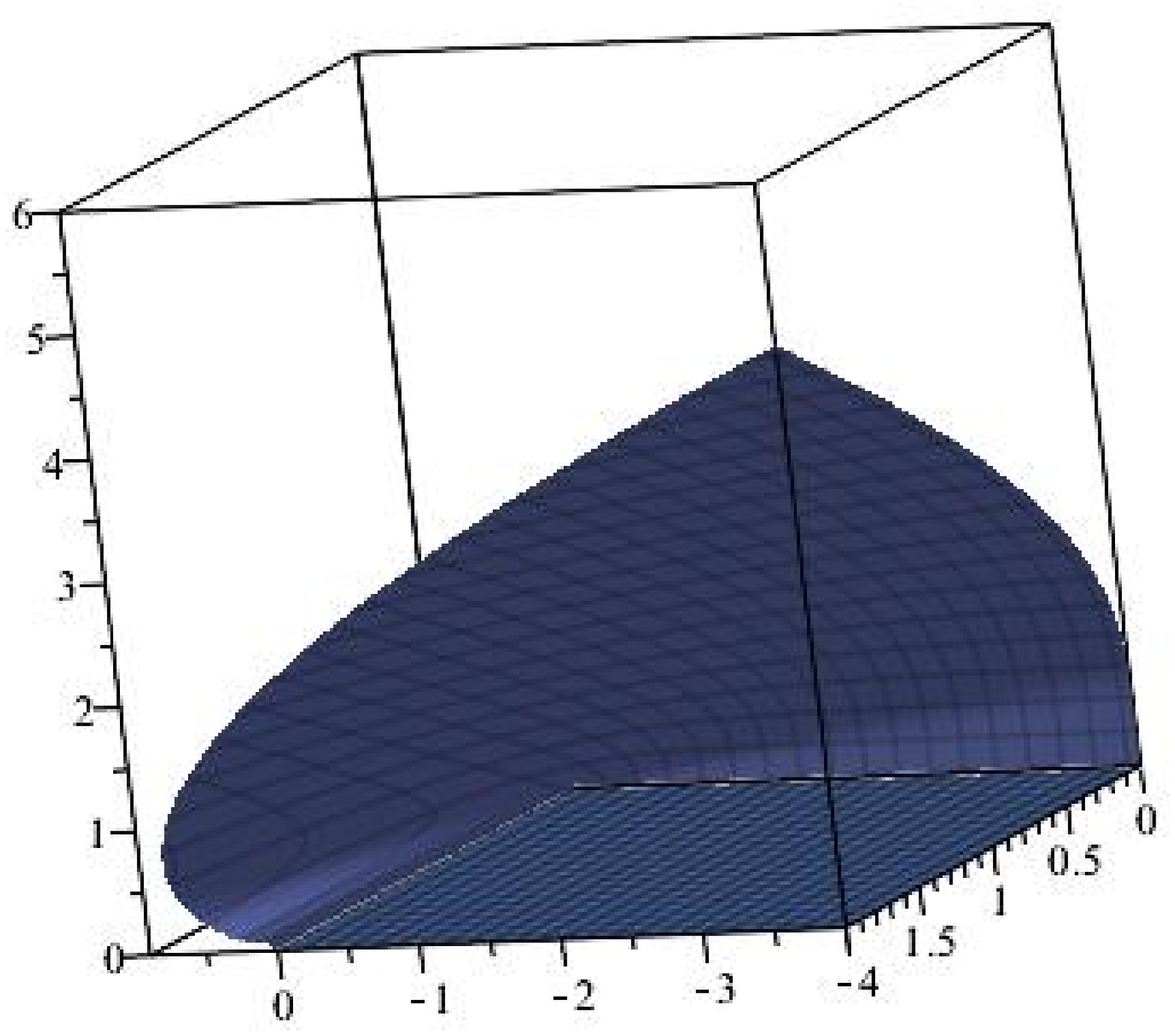}
\caption*{case 1}
\end{minipage}
\begin{minipage}[t]{0.3\textwidth}
\centering
\includegraphics[width=1.0 \columnwidth]{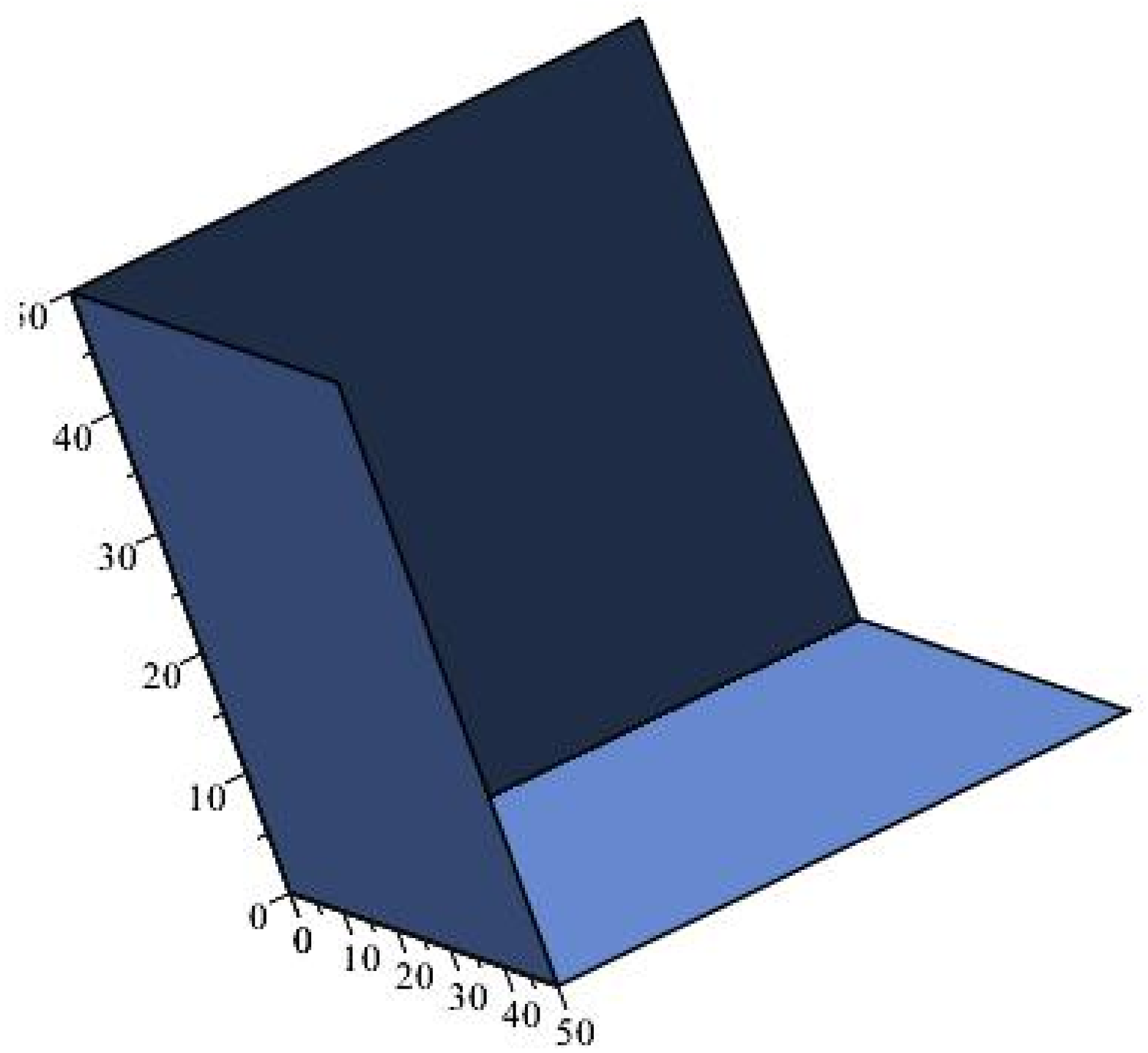}
\caption*{case 2}
\end{minipage}
\begin{minipage}[t]{0.3\textwidth}
\centering
\includegraphics[width=1.0 \columnwidth]{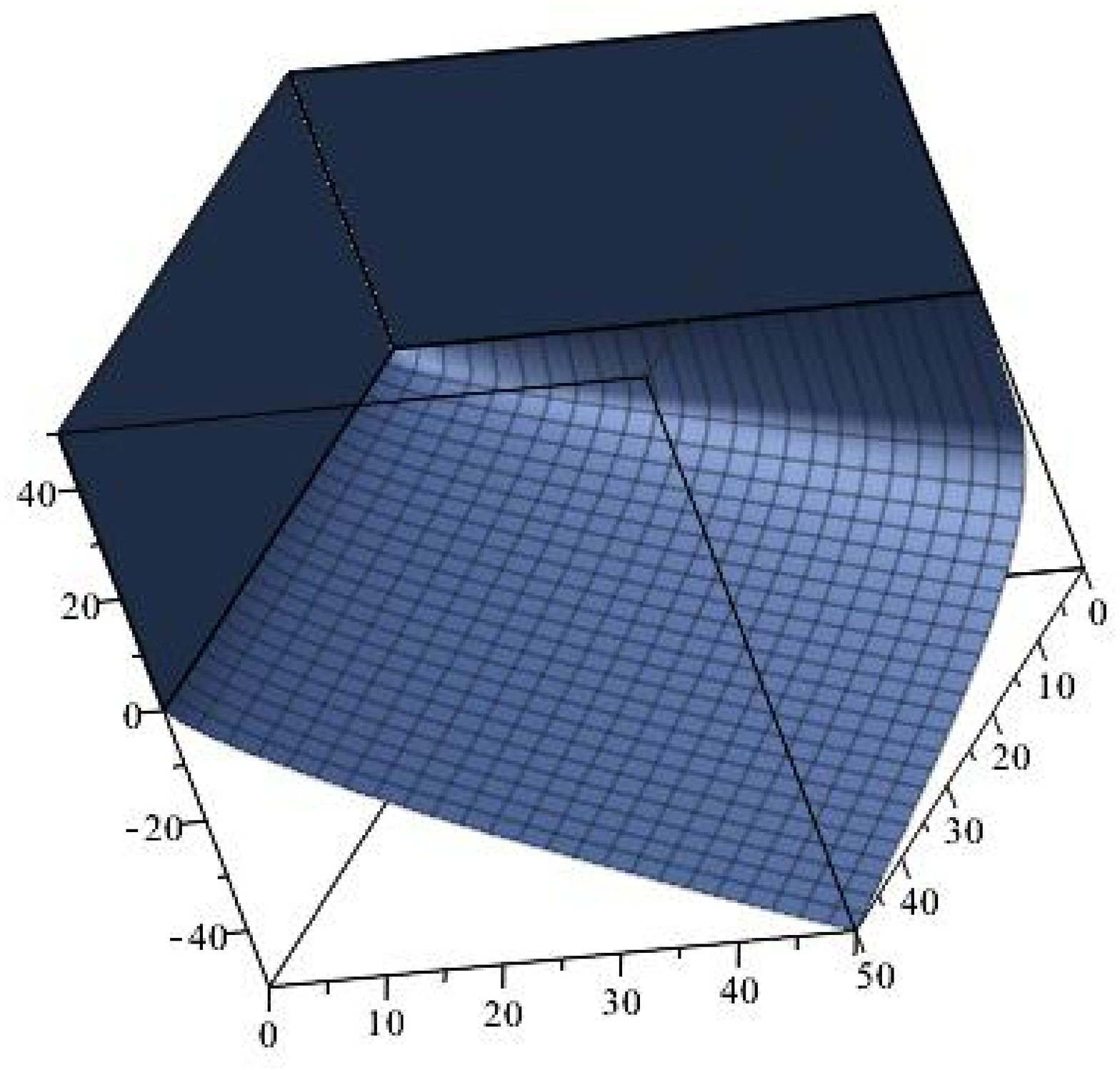}
\caption*{case 3}
\end{minipage}
}
\subfigure{
\begin{minipage}[t]{0.3\textwidth}
\centering
\includegraphics[width=1.0 \columnwidth]{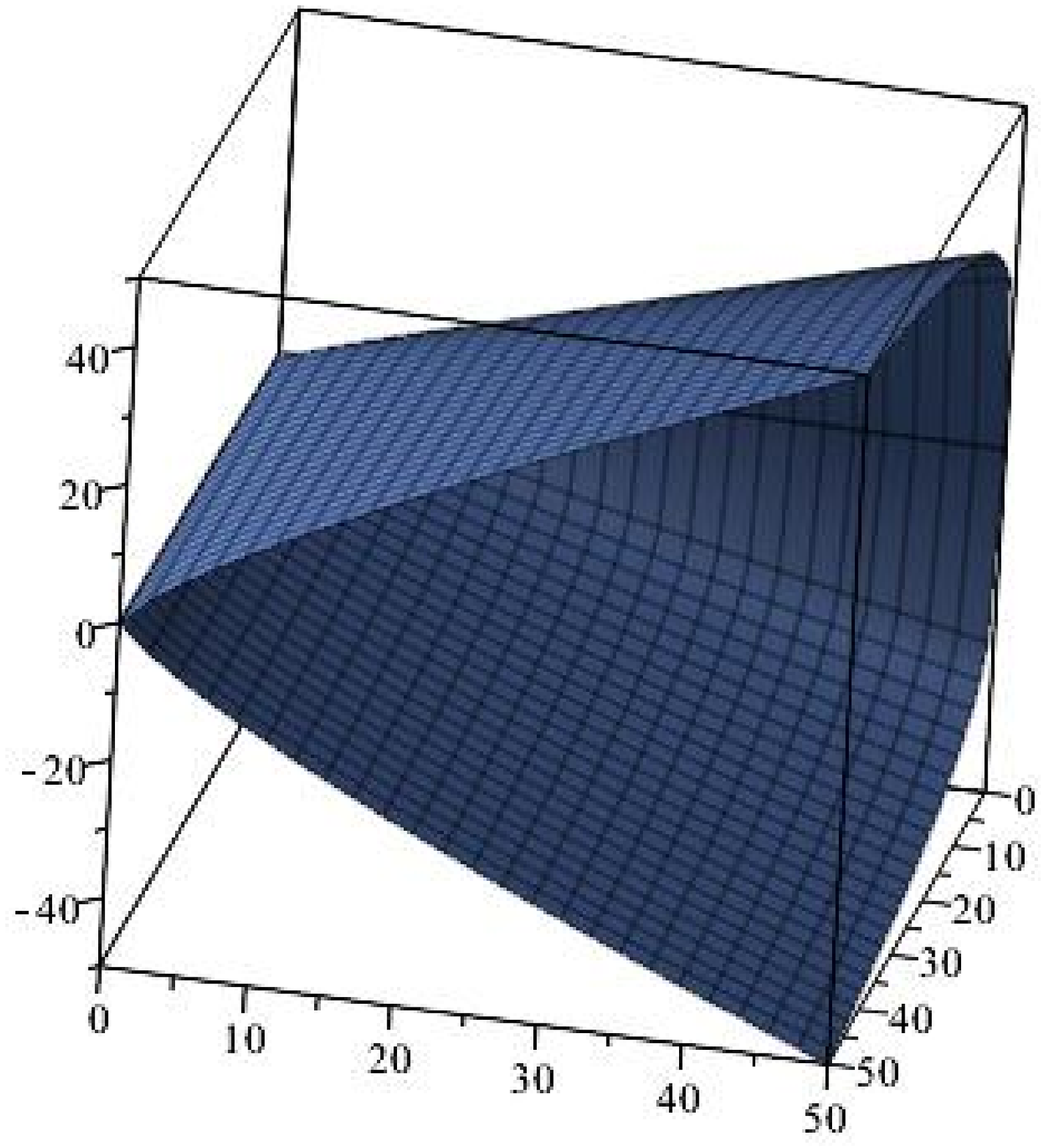}
\caption*{case 4}
\end{minipage}
\begin{minipage}[t]{0.3\textwidth}
\centering
\includegraphics[width=1.0 \columnwidth]{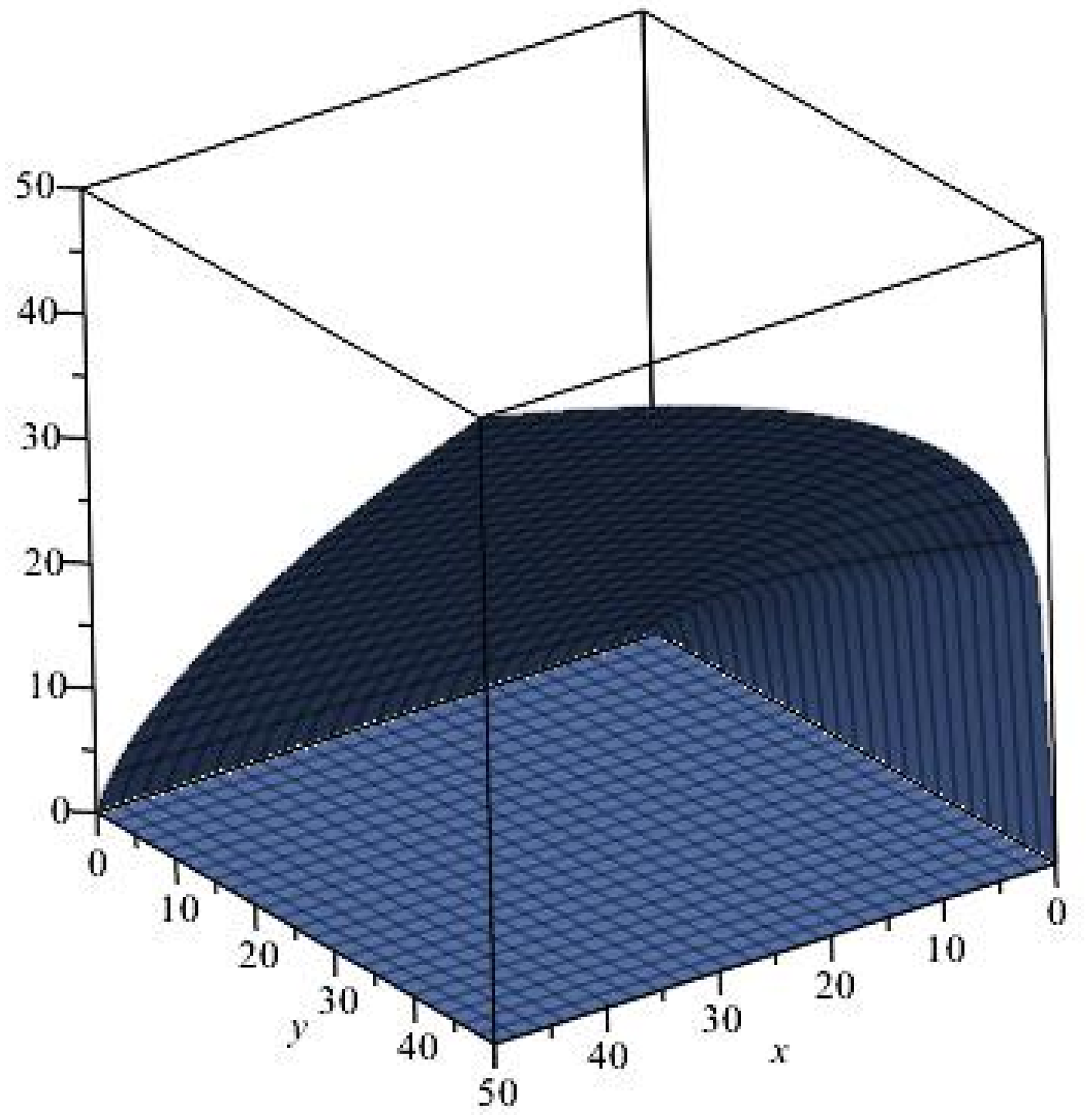}
\caption*{case 5}
\end{minipage}
}
\caption{Semi-homogeneous cone}
\end{figure}

\section{The associated families of affine spheres}
For convex affine spheres, the affine metric is positive-definite, and thus it provides a conformal structure in dimension two. There exists local complex coordinate  $z=x+i~y$ so that $g = 2 e^{\psi} ~dz d \bar z$. The affine cubic form $C = Udz^3+\bar U d\bar z^3$, where
$U=i|r_z, r_{zz},r|$. The affine normal satisfies:
\[
\xi = \Delta^g x/2 = e^{-\psi} x_{z \bar z}= -Hx.
\]
Simon-Wang derived in \cite{Simon-Wang.1993} the following equations for the frame $F=(r_{z},r_{\bar z},\xi)$:
\begin{align} \label{eqframe}
F_{z}=F\bpm  \psi_z & 0 & -H \\ -U e^{-\psi} & 0  & 0 \\ 0 & e^{\psi} & 0 \epm, ~~~
F_{\bar z}=F\bpm 0 & -\bar U e^{-\psi} & 0 \\ 0 &  \psi_{\bar z}  & -H \\ e^{\psi} & 0 & 0 \epm ;
\end{align}
and showed that the compatibility conditions (Tzitz\'{e}ica equation) are:
 \beq
\bca \label{eqtzi}
\psi_{z\bar z} +H e^{\psi} +|U|^2 e^{-2\psi} =0, & \\
U_{\bar z} =0 & .\\
\eca
 \eeq
The above system is invariant under $U \rightarrow e^{3it} U$ with $t \in [0, \frac{2}{3}\pi]$. Inserting the parameter $\l = e^{3it}$ into the frame equations \eqref{eqframe}, we obtain the following Lax representation of the system \eqref{eqtzi}:
\begin{align} \label{eqframe2}
F_{z}=F\bpm  \psi_z & 0 & -H \\  -\l U e^{-\psi} & 0  & 0 \\ 0 & e^{\psi} & 0 \epm, ~~~
F_{\bar z}=F\bpm 0 & -\l^{-1} \bar U e^{-\psi} & 0 \\ 0 &  \psi_{\bar z}  & -H \\ e^{\psi} & 0 & 0 \epm .
\end{align}

In this section we compute the isothermal parametrizations of Hildebrand's examples, then obtain the affine invariants. Substituting these affine invariants into the frame equations \eqref{eqframe2} and solving them, we get the whole associated family for each examples.

\subsection{The isothermal parametrization and the associated family of case 1}
In this subsection we consider the complete hyperbolic affine sphere asymptotic to the boundary of the cone obtained by the homogenization of the epigraph of the exponential function (case 1 in Theorem \ref{thm}). The corresponding example in \cite{Roland Hildebrand.{2013}} is given by:
\begin{align} \label{eqchas1}
 \mathbb{R}^2_{++} \ni \bpm x \\  y \epm \rightarrow  \vec{r} = \bpm r_1 \\  r_2 \\ r_3 \epm
 = \bpm -y(\frac{3}{2}\ln(x)+3\ln(y)+x) \\x^{-\frac{3}{2}}y^{-2}\sqrt{1+x} \\ y \epm,
\end{align}
here $\mathbb{R}_{++}$ is the set of positive reals.

Let us briefly review the classical results of isothermal coordinates in \cite{Spivak}.
If the metric is given locally as
\[
g = Edx^2+2Fdxdy+Gdy^2,
\]
then in the complex coordinate $z=x+i \ y$, it takes the form:
\[
g = \l \ |dz+\mu \ d \bar z|^2,
\]
with
\[
\l =\frac{1}{4}( E+G+2\sqrt{EG-F^2}), \;\;  \mu = \frac{E-G+2 i F}{4 \l}.
\]
In isothermal coordinates $(u, v)$ the metric should take the form
\[
g = \rho (d u^2+d v^2).
\]
The complex coordinate $\omega = u+i \ v$ satisfies
\[
\rho |d \omega|^2 = \rho |\omega_z|^2 |dz+\frac{\omega_{\bar z}}{\omega_z}d \bar z|^2,
\]
so that the coordinates $(u,v)$ will be isothermal if the \textbf{Beltrami equation}:
\begin{equation} \label{eqbeltrami}
\frac{\partial \omega}{\partial \bar z} = \mu \frac{\partial \omega}{\partial z}
\end{equation}
has a diffeomorphic solution.

To compute the isothermal parametrization of the affine sphere, we should compute the affine metric. Then we get
\[
L= \frac 18 \frac{(4x+3) (2x+3)^2}{x^{9/2} y (1+x)^{3/2}}, ~
M= \frac 34 \frac{(2x+3)^2}{x^{7/2} y^2 \sqrt{1+x}},~
N= \frac 32 \frac{(2x+3)^2}{x^{5/2} y^3 \sqrt{1+x}}.
\]
and
\[
\mu = \frac{\frac{\sqrt{3}}{6 \sqrt{x}\sqrt{1+x}}-\frac{1}{y}+\frac{i}{2x}}{\frac{\sqrt{3}}{6 \sqrt{x}\sqrt{1+x}}+\frac{1}{y}+\frac{i}{2x}}.
\]
The Beltrami equation \eqref{eqbeltrami} can be written as:
\begin{align} \label{eqbeltrami1}
d \omega &= \frac{\partial \omega}{\partial z} dz + \frac{\partial \omega}{\partial \bar z} d \bar z
         = (dz+\mu d \bar z) \frac{\partial \omega}{\partial z}
         = ((1+\mu) dx + i(1-\mu)dy) \frac{\partial \omega}{\partial z} \notag \\
         &= (\frac{2 \cdot (\frac{\sqrt{3}}{6 \sqrt{x}\sqrt{1+x}}+\frac{i}{2x})}{\frac{\sqrt{3}}{6 \sqrt{x}\sqrt{1+x}}+\frac{1}{y}+\frac{i}{2x}}dx
         +i \frac{2 \cdot (\frac{1}{y})}{\frac{\sqrt{3}}{6 \sqrt{x}\sqrt{1+x}}+\frac{1}{y}+\frac{i}{2x}}dy)
         \frac{\partial \omega}{\partial z}.
\end{align}
Then we divide our computation into the following steps: \\
(i) Set $d \tilde{x} = \frac{1}{2x}dx,\  d \tilde{y} = \frac{1}{y}dy$, i.e.  $\tilde{x} = \frac 12\ln(x), \ \tilde{y} = \ln (y)$. Then  \eqref{eqbeltrami1} becomes:
\begin{align} \label{eqbeltrami2}
d \omega = (\frac{2 \cdot (\frac{2 x \sqrt{3}}{6 \sqrt{x}\sqrt{1+x}}+i)}{\frac{\sqrt{3}}{6 \sqrt{x}\sqrt{1+x}}+\frac{1}{y}+\frac{i}{2x}}d\tilde{x}
         +i \frac{2}{\frac{\sqrt{3}}{6 \sqrt{x}\sqrt{1+x}}+\frac{1}{y}+\frac{i}{2x}}d\tilde{y})
         \frac{\partial \omega}{\partial z}.
\end{align}
and
\[
L=  \frac{(4e^{2\tilde{x}}+3) (2e^{2\tilde{x}}+3)^2}{e^{3\tilde{x}} (1+e^{2\tilde{x}})^{3/2}}, ~
M=
N=  \frac{3 (2e^{2\tilde{x}}+3)^2}{e^{3\tilde{x}} \sqrt{1+e^{2\tilde{x}}}}.
\]
(ii) Set $v=\tilde{y} - \tilde{x}$, we get
\begin{align} \label{eqbeltrami2}
d \omega = (\frac{2 \cdot (\frac{2 x \sqrt{3}}{6 \sqrt{x}\sqrt{1+x}})}{\frac{\sqrt{3}}{6 \sqrt{x}\sqrt{1+x}}+\frac{1}{y}+\frac{i}{2x}}d\tilde{x}
         +i \frac{2}{\frac{\sqrt{3}}{6 \sqrt{x}\sqrt{1+x}}+\frac{1}{y}+\frac{i}{2x}}dv)
         \frac{\partial \omega}{\partial z}.
\end{align}
\[
L=  \frac{(2e^{2\tilde{x}}+3)^2}{e^{\tilde{x}} (1+e^{2\tilde{x}})^{3/2}}, ~
M= 0, ~
N=  \frac{3 (2e^{2\tilde{x}}+3)^2}{e^{3\tilde{x}} \sqrt{1+e^{2\tilde{x}}}}.
\]
(iii) Set $d u=\frac{2 x \sqrt{3}}{6 \sqrt{x}\sqrt{1+x}} d\tilde{x} = \frac{\sqrt{3} e^{\tilde{x}}}{3 \sqrt{1+e^{2 \tilde{x}}}} d\tilde{x}$, i.e. $u= \frac{\text{arcsinh}(e^{\tilde{x}})}{\sqrt{3}}$. Then
\[
d \omega = du + i \ dv,
\]
and
\[
L = N = \frac{3\sqrt{3} (2\cosh^2(\sqrt{3}u)+1)^2}{(\cosh^2(\sqrt{3}u)-1)^2}, ~
M = 0.
\]
Summarize the above computation, we try the coordinate transformation:
\[
 x \rightarrow \sinh^2(\sqrt{3} x), \;\;y \rightarrow \sinh^{-1}( \sqrt{3} x)e^{y},
\]
and scale the immersion: $\vec{r} \rightarrow -\frac{1}{\sqrt{3}} \vec{r}$. Then
we can get the isothermal parametrization of the affine sphere:
\begin{align}
\vec{r'} = \bpm
\frac{1}{\sqrt{3}}\sinh^{-1}(\sqrt{3} x)e^{y} (\sinh^2(\sqrt{3} x)+3y) \\
 -\frac{1}{\sqrt{3}}\sinh^{-1}(\sqrt{3} x) e^{-2y} \sqrt{1+\sinh^2(\sqrt{3} x)}\\
   -\frac{1}{\sqrt{3}}\sinh^{-1}(\sqrt{3} x) e^{y}
   \epm.
\end{align}
By direct computation, we can get the affine metric, the affine mean curvature and the affine cubic form:
\[
g =3 (\csch(\sqrt{3}x)^2+\frac{2}{3})(dz~ d\bar z),~~~ H = -1, ~~~ U dz^3= i dz^3.
\]
It is easy to see that $g$ gives an explicit solution of Tzitz\'{e}ica equation \eqref{eqtzi}.
Substituting $g$, $H$, $U$ into frame equations \eqref{eqframe2} and solving them, then we can conclude the above computations into the following theorem.
\bthm
 The associated family of the complete hyperbolic affine sphere \eqref{eqchas1} in \cite{Roland Hildebrand.{2013}} is given by:
 \begin{equation}
 \vec{r}_{\l} = A \cdot (D\cdot B)^{-1} \cdot
 \begin{aligned}
 &\bpm \label{eqlambda}
 -(\coth(\sqrt{3}x)+\frac{2}{\sqrt{3}}\cos(t))e^{2(-\cos(t)x+\sin(t)y)} \\
 -(\coth(\sqrt{3}x)-\frac{-\sqrt{3}\sin(t)+\cos(t)}{\sqrt{3}})e^{(-\sqrt{3}\sin(t)+\cos(t))x+(-\sin(t)-\sqrt{3}\cos(t))y}\\
 -(\coth(\sqrt{3}x)-\frac{\sqrt{3}\sin(t)+\cos(t)}{\sqrt{3}})e^{(\sqrt{3}\sin(t)+\cos(t))x+(-\sin(t)+\sqrt{3}\cos(t))y}
  \epm, \\
 \end{aligned}
 \end{equation}
 here
  \begin{small}
 \[
 A = \left(
       \begin{array}{ccc}
         \frac{1}{\sqrt{3}} & \sqrt{2} & \frac{4}{\sqrt{3}} \ms \\
         -\frac{\sqrt{2}}{\sqrt{3}} & 1 & \frac{2\sqrt{2}}{\sqrt{3}} \ms \\
         - \frac{1}{\sqrt{3}} & \sqrt{2} & - \frac{1}{\sqrt{3}} \\
       \end{array}
     \right),
 D = \setlength\arraycolsep{0.3pt}
 \left(
       \begin{array}{ccc}
         (1\! +\sqrt{2})^{-\frac{2}{\sqrt{3}}\cos(t)} & 0 & 0 \ms \\
         0 & (1\! +\sqrt{2})^{-\sin(t)\! +\frac{1}{\sqrt{3}}\cos(t)} & 0 \ms \\
         0 & 0 & (1\! +\sqrt{2})^{\sin(t)\! +\frac{1}{\sqrt{3}}\cos(t)}\\
       \end{array}
     \right),
 \]

 \[
 \tiny{
 B =
 \left(
       \begin{array}{ccc}
         \! - \!(\sqrt{2}\! + \!\frac{2}{\sqrt{3}}\cos(t)) & (\! - \!\sqrt{3}\! + \!(\sqrt{2}\! + \!\frac{2}{\sqrt{3}}\cos(t))) \! \cdot \! (2\cos(t)) & \! - \!(\sqrt{2}\! + \!\frac{2}{\sqrt{3}}\cos(t)) \! \cdot \! (2\sin(t)) \ms \\
         \! - \!(\sqrt{2}\! + \!\sin(t)\! - \!\frac{\cos(t)}{\sqrt{3}}) & (\! - \!\sqrt{3}\! + \!(\sqrt{2}\! + \!\sin(t)\! - \!\frac{\cos(t)}{\sqrt{3}}) \! \cdot \! (\sqrt{3}\sin(t)\! - \!\cos(t)) & (\sqrt{2}\! + \!\sin(t)\! - \!\frac{\cos(t)}{\sqrt{3}}) \! \cdot \! (\sin(t))\! + \!\sqrt{3}\cos(t)) \ms \\
         \! - \!(\sqrt{2}\! - \!\sin(t)\! - \!\frac{\cos(t)}{\sqrt{3}}) &\! - \! (\! - \!\sqrt{3}\! + \!(\sqrt{2}\! - \!\sin(t)\! - \!\frac{\cos(t)}{\sqrt{3}}) \! \cdot \! (\sqrt{3}\sin(t)\! + \!\cos(t)) & (\sqrt{2}\! - \!\sin(t)\! - \!\frac{\cos(t)}{\sqrt{3}}) \! \cdot \! (\sin(t))\! - \!\sqrt{3}\cos(t))\\
       \end{array}
     \right),
      }
 \]
\end{small}
with the affine metric, the affine mean curvature and the affine cubic form:
\[
g =3 (\csch(\sqrt{3}x)^2+\frac{2}{3})(dz~ d\bar z),~~~ H = -1, ~~~ U dz^3 = e^{3it} dz^3.
\]
\ethm
 \brem
 The expression \eqref{eqlambda} is gauged by the initial condition
 \[
 (\vec{r_{\l}}, (\vec{r_{\l}})_x, (\vec{r_{\l}})_y)|_{x=\frac{\ln(1+\sqrt{2})}{\sqrt{3}}, \ y=0} := (\vec{r'}, \vec{r'}_x, \vec{r'}_y)|_{x=\frac{\ln(1+\sqrt{2})}{\sqrt{3}}, \ y=0} = A.
 \]
 With this initial condition, one can check $\lim\limits_{t \to \frac{\pi}{6}}\vec{r_{\l}} = \vec{r'}$ by direct computation.
 \erem
 \subsection{The isothermal parametrizations of remaining cases}
The complete hyperbolic affine spheres asymptotic to the boundary of the remaining cones can be treated in a common framework. In this subsection we compute the common isothermal coordinates transformation for these affine spheres by Weierstrass elliptic functions. Let us review  the results in \cite{Roland Hildebrand.{2013}}.

Let $p \in [2,\infty)$, $q \in (1,2]$  be reals such that $\frac{1}{p}+\frac{1}{q}=1$, and let $\a, ~\beta \in [ 0, +\infty ]$. The cases $2-5$ in Theorem \ref{thm} are given by the closure of
\[
K^\circ = \{(x,y,z)^T | -\a x^{1/p} y^{1/q} < z < \beta x^{1/p} y^{1/q}, x>0,y>0\},
\]
with the value $(\a,\beta) = (0, \infty), (\infty, 1), (\a,1), (0,1)$, respectively. The PDE \eqref{eqMA2} becomes the ODE
\begin{align}
e^{-\phi} \dot{\phi} (t \dot{\phi}+p+1) (t \dot{\phi}+q+1) = (p+q)te^{\phi}+c ,\label{eqode}\\
t \dot{\phi} > -\frac{2}{3}(p+q)-\frac{1}{3} ~, ~~~ \ddot{\phi} > \frac{\dot{\phi}^2(p+q-1+t\dot{\phi})}{2(p+q)+1+3t\dot{\phi}} ~,
\end{align}
where $\phi(t): (-\a, \beta)\rightarrow \mathbb{R}$ is a function of variable $t$, and $c$ is an integration constant. Suppose $t(\xi)$ is a function of $\xi$, with $\xi$ is a real parameter taking values in some interval $\Xi \subset \mathbb{R}$. Then the affine sphere is given by the immersion
 \begin{align} \label{eqchasg1}
 \Xi \times \mathbb{R} \ni \bpm \xi \\  \mu \epm \rightarrow  \vec{r} = \bpm r_1 \\  r_2 \\ r_3 \epm
 =e^{\frac{\phi(\xi)}{3}}
 \bpm  e^{\frac{q+1}{3q}\mu} \\
   e^{-\frac{p+1}{3p}\mu} \\
   e^{-\frac{p-q}{3(p+q)}\mu} t(\xi)\epm.
 \end{align}

Introduce the variables $\tau = \ln(t)$, $\varphi = \phi+\tau$, $\xi = \frac{d\varphi}{d\tau}=1+t\dot{\phi}$ and denote $P = (\xi-1) (\xi+p) (\xi+q)$ for writing shorthand. By virtue of \eqref{eqode}, we obtain
 \begin{align}\label{eqchasg2}
 e^\varphi = \frac{-c+s \sqrt{c^2+4 (p+q)P}}{2(p+q)},
 \end{align}
 \begin{align}\label{eqchasg3}
\frac{d \tau}{d\xi} = \frac{2(p+q)(3\xi+2(p+q-1))}{\sqrt{c^2+4 (p+q)P}(-c s+\sqrt{c^2+4 (p+q)P}}
  \end{align}
where  $s = \pm 1$. It is easy to see that if $\phi(t)$ is a solution of \eqref{eqode} for the value c, then $\tilde{\phi}(t) = \phi(-t)$ is a solution for $-c$. We may thus assume without restriction of generality that $c \leq 0$ and have the following classification: \\
(1) when $\xi \equiv c_1$ is constant, the affine spheres asymptotic to the orthant (case 2 of Theorem 2.1) are the level surfaces of the product $xyz$. \\
(2) when $c=-2(p+q)$, $s = -1$, the affine spheres asymptotic to the cones given in the case 3 of Theorem 2.1, up to a sign change in the second coordinate.\\
(3) when $c=-2(p+q)$, $s = 1$, the affine spheres asymptotic to the cones given in the case 5 of Theorem 2.1. \\
(4) when $c=0,s=1$, the affine spheres asymptotic to the cones given in the case 4 of Theorem 2.1 with $\a = 1$. \\
(5) when $-2(p+q)<c<0$, the affine spheres asymptotic to the cones given in the case 4 of Theorem 2.1 for $\a \neq 1$.

Now we begin to compute the transformation. We consider Blaschke metric of the immersion
\begin{align*}
g = \frac{Ldx^2+2Mdxdy+Ndy^2}{|LN-M^2|^{\frac{1}{4}}},
\frac{}{}\end{align*}
here
\begin{align*}
& L = |r_\xi,r_\mu,r_{\xi \xi}|=e^\varphi [-\frac{1}{9} (\frac{d\tau}{d\xi})^3 (\xi^2+\xi-2)+\frac{1}{3}(\frac{d\tau}{d\xi})^2], \\
& M = |r_\xi,r_\mu,r_{\xi \mu}|=\frac{p-2}{9p} e^\varphi (\frac{d\tau}{d\xi})^2 (\xi-1), \\
& N = |r_\xi,r_\mu,r_{\mu \mu}|=\frac{p-1}{9p^2} e^\varphi (\frac{d\tau}{d\xi}) (3\xi-2(p+q-1)).
\end{align*}
Similar to case 1, we divide the computation into the following 3 steps: \\
(i) Set $d\xi = \frac{3\xi+2(p+q-1)}{2(\xi-1)}\frac{d\xi}{d\tau} d\xi_1$, $\mu = \frac{p(p-2)}{2(p-1)}y_1$, the coefficients of metric become:
\begin{align*}
& L = -\frac{p(p-2)}{36(p-1)}e^\varphi \frac{(3\xi+2(p+q-1))^3}{(\xi-1)^3}[(\xi^2+\xi-2)-3(\frac{d\xi}{d\tau})] \\
& M = N = \frac{p(p-2)^3}{36(p-1)^2}e^\varphi \frac{(3\xi+2(p+q-1))^2}{\xi-1}
\end{align*}
(ii) Set $y_1=y_2-\xi_1$, we get
\begin{align*}
& L = -\frac{p(p-2)}{12(p-1)}e^\varphi \frac{(3\xi+2(p+q-1))^2}{(\xi-1)^3}(\sqrt{\frac{c^2}{4(p+q)}+P}-\frac{c s }{\sqrt{4(p+q)}})^2 \\
& M = 0, \;\; N = \frac{p(p-2)^3}{36(p-1)^2}e^\varphi \frac{(3\xi+2(p+q-1))^2}{\xi-1}
\end{align*}
(iii) Set $d\xi_1 = \frac{\xi-1}{\sqrt{\frac{c^2}{4(p+q)}+P}-\frac{c s}{\sqrt{4(p+q)}}} d\xi_2$, $y_2=\frac{\sqrt{3(p-1)}}{p-2}y_3$, we finally arrive at
\begin{align*}
L = N = \frac{p\sqrt{3(p-1)}}{12(p-1)} e^\varphi \frac{(3\xi+2(p+q-1))^2}{\sqrt{\frac{c^2}{4(p+q)}+P}-\frac{c s}{\sqrt{4(p+q)}}}; \;\;\; M=0.
\end{align*}
Summarizing the above computation, we obtain the transformation:
\begin{align}
d\xi = \sqrt{\frac{c^2}{4(p+q)}+P} d\xi_2, \;\;\;  \mu = \frac{p(p-2)}{2(p-1)}(\frac{\sqrt{3(p-1)}}{p-2}y_3-\xi_1).
\end{align}
The above integral can be expressed by Weierstrass elliptic function:
\begin{align}
\xi = 4 \mathscr{P} (\xi_2+\omega_1,g_2,g_3)-b_1,
\end{align}
here we use $\frac{c^2}{4(p+q)}+P = \xi^3+3b_1\xi^2+3b_2\xi+b_3$, $g_2=\frac{3}{4}(b_1^2-b_2)$, $g_3=\frac{1}{16}(3b_1b_2-2b_1^3-b_3)$. $\mathscr{P}(\xi_2,g_2,g_3)$ is Weierstrass P-function with the invariants $g_2$, $g_3$, and $\omega_1$ is the first half-period of $\mathscr{P}$. We use $\mathscr{P}(\xi_2)$ to denote $\mathscr{P}(\xi_2,g_2,g_3)$ for convenient. This transformation allows to compute the isothermal parametrizations of the affine spheres:
\begin{eqnarray*}
r_1 &=& e^{\frac{\varphi(\xi)}{3}} t^{-\frac{1}{3}} e^{\frac{q+1}{3q}\mu} =  e^{\frac{\varphi(\xi)-\tau}{3}} e^{\frac{q+1}{3q}\mu} \notag\\
    &=& e^{\frac{\varphi(\xi)-\tau}{3}-\frac{q+1}{3q}\frac{p(p-2)}{2(p-1)}\xi_1} e^{\frac{q+1}{3q}\frac{p(p-2)}{2(p-1)}\frac{\sqrt{3(p-1)}}{p-2}y_3} \notag \\
    &=&e^{\int_0^{\xi_2}(\frac{1}{3}(\frac{d\varphi(\xi)}{d\xi_2}-\frac{d\tau}{d\xi_2})
    -\frac{q+1}{3q}\frac{p(p-2)}{2(p-1)}\frac{d\xi_1}{d\xi_2})d\xi_2}
    e^{\frac{(2p-1)}{2\sqrt{3(p-1)}}y_3} \notag\\
    &=& e^{\frac{1}{2}\int_0^{\xi_2}\frac{(\xi-1)(\xi+q)}{\sqrt{\frac{c^2}{4(p+q)}+P}-\frac{c s}{\sqrt{4(p+q)}}}d\xi_2} e^{\frac{(2p-1)}{2\sqrt{3(p-1)}}y_3} \notag \\
    &=& \sqrt{\frac{4 \mathscr{P} (\xi_2+\omega_1)-b_1+p}{4 \mathscr{P} (\omega_1)-b_1+p}
        \left(\frac{\sigma(\xi_2+\omega_1+a_1)\sigma(\omega_1-a_1)}{\sigma(\xi_2+\omega_1-a_1)\sigma(\omega_1+a_1)}\right)^{-s}}
        e^{s\zeta(a_1)\xi_2+\frac{(2p-1)}{2\sqrt{3(p-1)}}y_3}
\end{eqnarray*}
here $a_1$ is the root of $4 \mathscr{P} (\xi_2,g_2,g_3)-b_1+p=0$, and $\sigma$, $\zeta$ are Weierstrass sigma function and Weierstrass zeta function respectively.
By similar computation, we have
\begin{eqnarray}
r_1 = \sqrt{\frac{4 \mathscr{P} (\xi_2+\omega_1)-b_1+p}{4 \mathscr{P} (\omega_1)-b_1+p}
        \left(\frac{\sigma(\xi_2+\omega_1+a_1)\sigma(\omega_1-a_1)}{\sigma(\xi_2+\omega_1-a_1)\sigma(\omega_1+a_1)}\right)^{-s}}
        e^{s\zeta(a_1)\xi_2+\frac{(2p-1)}{2\sqrt{3(p-1)}}y_3} \label{eqas1} \\
r_2 = \sqrt{\frac{4 \mathscr{P} (\xi_2+\omega_1)-b_1+q}{4 \mathscr{P} (\omega_1)-b_1+q}
        \left(\frac{\sigma(\xi_2+\omega_1+a_2)\sigma(\omega_1-a_2)}{\sigma(\xi_2+\omega_1-a_2)\sigma(\omega_1+a_2)}\right)^{-s}}
        e^{s\zeta(a_2)\xi_2-\frac{(p+1)}{2\sqrt{3(p-1)}}y_3} \label{eqas2}\\
r_3 = \sqrt{\frac{4 \mathscr{P} (\xi_2+\omega_1)-b_1-1}{4 \mathscr{P} (\omega_1)-b_1-1}
        \left(\frac{\sigma(\xi_2+\omega_1+a_3)\sigma(\omega_1-a_3)}{\sigma(\xi_2+\omega_1-a_3)\sigma(\omega_1+a_3)}\right)^{-s}}
        e^{s\zeta(a_3)\xi_2-\frac{(p-2)}{2\sqrt{3(p-1)}}y_3} \label{eqas3}
\end{eqnarray}
here $a_2$, $a_3$ are the roots of $4 \mathscr{P} (\xi_2,g_2,g_3)-b_1+q=0$ and $4 \mathscr{P} (\xi_2,g_2,g_3)-b_1-1=0$ respectively.  It is easy to check that the Weierstrass P-function can be expressed by Weierstrass sigma function:
\begin{equation}
\mathscr{P}(\xi_2+\omega_1)-\mathscr{P}(a_i) = -\frac{\sigma(\xi_2+\omega_1+a_i)\sigma(\xi_2+\omega_1-a_i)}{\sigma(\xi_2+\omega_1)^2\sigma(a_i)^2}.
\end{equation}
By this expression and a gauge transformation:
\[
\bpm r_1 \\  r_2 \\ r_3 \epm  \rightarrow  \left(
                                             \begin{array}{ccc}
                                               \frac{1}{3p} & 0 & 0 \\
                                               0 & \frac{1}{3p} & 0 \\
                                               0 & 0 & \frac{\sqrt{3}(p-1)c}{2} \\
                                             \end{array}
                                           \right)
\cdot \bpm r_1 \\  r_2 \\ r_3 \epm,
\]
we can simplify \eqref{eqas1}-\eqref{eqas3} as
\begin{eqnarray}
r_1 &=& \frac{1}{3p}\left(\frac{\sigma(\xi_2+\omega_1-s \cdot a_1)\sigma(\omega_1)}{\sigma(\xi_2+\omega_1)\sigma(\omega_1-s \cdot a_1)}\right)
        e^{s\zeta(a_1)\xi_2+\frac{(2p-1)}{2\sqrt{3(p-1)}}y_3}, \label{eqaffsim1}\\
r_2 &=& \frac{1}{3p}\left(\frac{\sigma(\xi_2+\omega_1-s \cdot a_2)\sigma(\omega_1)}{\sigma(\xi_2+\omega_1)\sigma(\omega_1-s \cdot a_2)}\right)
        e^{s\zeta(a_2)\xi_2-\frac{(p+1)}{2\sqrt{3(p-1)}}y_3}, \label{eqaffsim2}\\
r_3 &=& \frac{\sqrt{3}(p-1)c}{2}\left(\frac{\sigma(\xi_2+\omega_1-s \cdot a_3)\sigma(\omega_1)}{\sigma(\xi_2+\omega_1)\sigma(\omega_1-s \cdot a_3)}\right)
        e^{s\zeta(a_3)\xi_2-\frac{(p-2)}{2\sqrt{3(p-1)}}y_3}. \label{eqaffsim3}
\end{eqnarray}
\brem
Since $\mathscr{P}$ is even, both of $a_i$ and $-a_i$ (i=1,2,3) are the roots of those  equations. We need to choose $a_i$ such that $\mathscr{P}' (a_i,g_2,g_3)<0$ .
\erem

Then we obtain the affine metric, the affine mean curvature and the affine cubic form:
\begin{eqnarray}
g = (\mathscr{P}(\xi_2+\omega_1)+\frac{2(p+q-1)-3b_1}{12}) (d\xi_2^2+dy_3^2), ~~ H = -1, \label{eqinvirant1}\\
 U d(\xi_2+i y_3)^3= s \cdot \frac{c \cdot \sqrt{p-1}}{32p}-i \frac{(p-2)(2p-1)(p+1)}{2(12(p-1))^{\frac{3}{2}}}d(\xi_2+i y_3)^3. \label{eqinvirant2}
\end{eqnarray}
The invariants \eqref{eqinvirant1}-\eqref{eqinvirant2} imply that the two affine spheres coresponding to $s=1$ and $s=-1$ are in the same associated family. The coefficients $U$ of these two affine spheres can be expressed as $|U|e^{3it}$ and $|U|e^{\pi-3it}$ respectively, in other words, $s = -Sgn(cos(3t))$.

\subsection{The associated families of remaining cases}
With the generic formulas \eqref{eqaffsim1}-\eqref{eqaffsim3}, we can construct the complete hyperbolic affine spheres and the associated families case by case.
Some parts of construction  have referred to the formulas of the affine spheres in \cite{Roland Hildebrand.{2013}}. We consider the following cases.

$\mathbf{-2(p+q)<c<0.}$ In this case $\beta=1$. For positive t, $s=1$, \eqref{eqaffsim1}-\eqref{eqaffsim3} yielding
\begin{align} \label{eqcase4g1}
\mathbb{R}_{--} \times \mathbb{R} \ni \bpm \xi_2 \\  y_3 \epm \rightarrow  \vec{r} =
\bpm
\frac{1}{3p}\left(\frac{\sigma(\xi_2+\omega_1-a_1)\sigma(\omega_1)}{\sigma(\xi_2+\omega_1)\sigma(\omega_1-a_1)}\right)
        e^{\zeta(a_1)\xi_2+\frac{(2p-1)}{2\sqrt{3(p-1)}}y_3} \\
\frac{1}{3p}\left(\frac{\sigma(\xi_2+\omega_1-a_2)\sigma(\omega_1)}{\sigma(\xi_2+\omega_1)\sigma(\omega_1-a_2)}\right)
        e^{\zeta(a_2)\xi_2-\frac{(p+1)}{2\sqrt{3(p-1)}}y_3}\\
\frac{\sqrt{3}(p-1)c}{2}\left(\frac{\sigma(\xi_2+\omega_1-a_3)\sigma(\omega_1)}{\sigma(\xi_2+\omega_1)\sigma(\omega_1-a_3)}\right)
        e^{\zeta(a_3)\xi_2-\frac{(p-2)}{2\sqrt{3(p-1)}}y_3}
\epm.
\end{align}
Since $\phi$ is analytic in $t$, it must be an even function of $t$. For negative t, $s=-1$, we have
\begin{align}\label{eqcase4g2}
\mathbb{R}_{--} \times \mathbb{R} \ni \bpm \xi_2 \\  y_3 \epm \rightarrow  \vec{r} =
\bpm
\frac{1}{3p}\left(\frac{\sigma(\xi_2+\omega_1+a_1)\sigma(\omega_1)}{\sigma(\xi_2+\omega_1)\sigma(\omega_1+a_1)}\right)
        e^{-\zeta(a_1)\xi_2+\frac{(2p-1)}{2\sqrt{3(p-1)}}y_3} \\
\frac{1}{3p}\left(\frac{\sigma(\xi_2+\omega_1+a_2)\sigma(\omega_1)}{\sigma(\xi_2+\omega_1)\sigma(\omega_1+a_2)}\right)
        e^{-\zeta(a_2)\xi_2-\frac{(p+1)}{2\sqrt{3(p-1)}}y_3}\\
\frac{\sqrt{3}(p-1)c}{2}\left(\frac{\sigma(\xi_2+\omega_1+a_3)\sigma(\omega_1)}{\sigma(\xi_2+\omega_1)\sigma(\omega_1+a_3)}\right)
        e^{-\zeta(a_3)\xi_2-\frac{(p-2)}{2\sqrt{3(p-1)}}y_3}
\epm.
\end{align}
This affine sphere is asymptotic to the cone given in case 4 of Theorem \ref{thm}, for $\a\neq1$. The affine invariants are given in \eqref{eqinvirant1}-\eqref{eqinvirant2}. When computing the associated family, we can consider the coefficient $U$ of this affine sphere is $|U|\cdot \arg(U)$. Then, for convenience, substitute $g$, $H$, $U'=|U|$ into frame equations \eqref{eqframe2} and solve them. We conclude the above discussions into the following theorem.
\bthm The associated family the complete hyperbolic affine sphere \eqref{eqchasg1}-\eqref{eqchasg3} with $-2 (p+q) < c <0$ in \cite{Roland Hildebrand.{2013}} is given by
\begin{align} \label{eqfamily1}
 \vec{r}_{\l} =
 \bpm \ms
\frac{1}{3p}\left(\frac{\sigma(\xi_2+\omega_1-s \cdot m_1)\sigma(\omega_1)}{\sigma(\xi_2+\omega_1)\sigma(\omega_1-s \cdot m_1)}\right)
        e^{s\zeta(m_1)\xi_2+f_1(t)y_3}\\ \ms
\frac{1}{3p}\left(\frac{\sigma(\xi_2+\omega_1-s \cdot m_2)\sigma(\omega_1)}{\sigma(\xi_2+\omega_1)\sigma(\omega_1-s \cdot m_2)}\right)
        e^{s\zeta(m_2)\xi_2+f_2(t)y_3}\\
\frac{\sqrt{3}(p-1)c}{2}\left(\frac{\sigma(\xi_2+\omega_1-s \cdot m_3)\sigma(\omega_1)}{\sigma(\xi_2+\omega_1)\sigma(\omega_1-s \cdot m_3)}\right)
        e^{s\zeta(m_3)\xi_2+f_3(t)y_3}
\epm,
\end{align}
where $s=-Sgn(\cos(3t))$ and $m_i$ are the solutions of
\begin{equation*}
\mathscr{P} (\xi_2)+\frac{b_1}{4}+\frac{-\sqrt{-2c(p-1)}U'\sin(3t)}{3^{3/4}p\cdot f_i(t)}=0,  ~~ \mathscr{P}'(m_i)<0.
\end{equation*}
The functions $f_i(t)$ are given by
\begin{align*}
f_1(t)& = -\frac{1}{\sqrt{3}} \cdot \sqrt{\frac{p^2-p+1}{p-1}} \sin(\frac{1}{6} \pi+ \frac{1}{3} \arccos(\frac{U'H\sin(3t)}{\sqrt{\left(\frac{p^2-p+1}{12\cdot (p-1)}\right)^3}}), \\
f_2(t)& = -\frac{1}{\sqrt{3}} \cdot \sqrt{\frac{p^2-p+1}{p-1}} \cos(\frac{1}{3} \pi+ \frac{1}{3} \arccos(\frac{U'H\sin(3t)}{\sqrt{\left(\frac{p^2-p+1}{12\cdot (p-1)}\right)^3}}), \\
f_3(t)& = \frac{1}{\sqrt{3}} \cdot \sqrt{\frac{p^2-p+1}{p-1}} \cos( \frac{1}{3} \arccos(\frac{U'H\sin(3t)}{\sqrt{\left(\frac{p^2-p+1}{12\cdot (p-1)}\right)^3}}).
\end{align*}
What's more, the affine metric, the affine mean curvature and the affine cubic form are given as
\begin{align*}
g = (\mathscr{P}(\xi_2+\omega_1)+\frac{2(p+q-1)-3b_1}{12}) (d\xi_2^2+dy_3^2),~~ H = -1, \\
U d(\xi_2+iy_3)^3= e^{3it}\cdot \sqrt{\frac{c^2\cdot (p-1)}{(32p)^2}+\frac{(p-2)^2(2p-1)^2(p+1)^2}{4(12(p-1))^3}} d(\xi_2+iy_3)^3.
\end{align*}
\ethm
$\mathbf{c=0.}$ In this case $s=1$, $\beta=1$ and \eqref{eqaffsim1}-\eqref{eqaffsim3} the become
\begin{align} \label{eqcase4s1}
\mathbb{R}_{--} \times \mathbb{R} \ni \bpm \xi_2 \\  y_3 \epm \rightarrow  \vec{r} = \bpm r_1 \\  r_2 \\ r_3 \epm
 = \tiny{\bpm
\frac{2}{3p\sqrt{p+1}}\sqrt{\mathscr{P}(\xi_2+\omega_1)+\frac{-b_1+p}{4}}e^{\frac{2p-1}{2\sqrt{3(p-1)}}y_3} \\
\frac{2}{3\sqrt{p(2p-1)}}\sqrt{\mathscr{P}(\xi_2+\omega_1)+\frac{-b_1+q}{4}}e^{-\frac{p+1}{2\sqrt{3(p-1)}}y_3} \\
2\sqrt{3p(p-1)(p+1)(2p-1)}\sqrt{\mathscr{P}(\xi_2+\omega_1)+\frac{-b_1-1}{4}}e^{-\frac{p-2}{2\sqrt{3(p-1)}}y_3}
\epm},
\end{align}
for positive t. Then for negative  t, we have
\begin{align} \label{eqcase4s2}
\mathbb{R}_{--} \times \mathbb{R} \ni \bpm \xi_2 \\  y_3 \epm \rightarrow  \vec{r} = \bpm r_1 \\  r_2 \\ r_3 \epm
 = \tiny{\bpm
\frac{2}{3p\sqrt{p+1}}\sqrt{\mathscr{P}(\xi_2+\omega_1)+\frac{-b_1+p}{4}}e^{\frac{2p-1}{2\sqrt{3(p-1)}}y_3} \\
\frac{2}{3\sqrt{p(2p-1)}}\sqrt{\mathscr{P}(\xi_2+\omega_1)+\frac{-b_1+q}{4}}e^{-\frac{p+1}{2\sqrt{3(p-1)}}y_3} \\
-2\sqrt{3p(p-1)(p+1)(2p-1)}\sqrt{\mathscr{P}(\xi_2+\omega_1)+\frac{-b_1-1}{4}}e^{-\frac{p-2}{2\sqrt{3(p-1)}}y_3}
\epm}.
\end{align}
This affine sphere is asymptotic to the cone given in case 4 of Theorem \ref{thm}, with $\a=1$. The affine metric, the affine mean curvature and the affine cubic form of this affine sphere are
\begin{eqnarray*}
g = (\mathscr{P}(\xi_2+\omega_1,g_2,g_3)+\frac{2(p+q-1)-3b_1}{12}) (d\xi_2^2+dy_3^2), \\
H = -1, \;\;~~ U d(\xi_2+i y_3)^3=-i \frac{(p-2)(2p-1)(p+1)}{2(12(p-1))^{\frac{3}{2}}}d(\xi_2+i y_3)^3.
\end{eqnarray*}

\brem
By direct computation, it is easy check the expression of affine spheres \eqref{eqcase4s1}-\eqref{eqcase4s2} can get from  \eqref{eqcase4g1}-\eqref{eqcase4g2} when $c \rightarrow 0$. Henceforth, the associated families in this case can also get from \eqref{eqfamily1} by  $c \rightarrow 0$.
\erem

$\mathbf{c=-2(p+q),\beta=1.}$ In this case $s=1$, and \eqref{eqaffsim1}-\eqref{eqaffsim3} the becomes

\begin{align*}
&\mathbb{R}_{--} \times \mathbb{R}  \ni \bpm \xi_2 \\  y_3 \epm \rightarrow  \vec{r} = \bpm r_1 \\  r_2 \\ r_3 \epm
 =\\
& \bpm
\frac{1}{\sqrt{3(p-1)}} (\sqrt{p^2-p+1} \coth(\frac{\sqrt{p^2-p+1}}{2\sqrt{p-1}}\xi_2)+1)e^{-\frac{1}{2\sqrt{p-1}}\xi_2+\frac{2p-1}{2\sqrt{3(p-1)}}y_3} \\
\frac{1}{\sqrt{3(p-1)}} (\sqrt{p^2-p+1} \coth(\frac{\sqrt{p^2-p+1}}{2\sqrt{p-1}}\xi_2)+(p-1))e^{-\frac{\sqrt{p-1}}{2}\xi_2-\frac{p+1}{2\sqrt{3(p-1)}}y_3} \\
-\frac{1}{\sqrt{3}p} (\sqrt{p^2-p+1} \coth(\frac{\sqrt{p^2-p+1}}{2\sqrt{p-1}}\xi_2)-p)e^{\frac{p}{2\sqrt{p-1}}\xi_2-\frac{p-2}{2\sqrt{3(p-1)}}y_3}
\epm.
\end{align*}
This affine sphere is asymptotic to the cone given in case 5 of Theorem \ref{thm}. The affine metric, the affine mean curvature and the affine cubic form of this affine sphere are
\begin{eqnarray*}
g = \frac {p^2-p+1}{4(p-1)}(\coth^2(\frac{\sqrt{p^2-p+1}}{2\sqrt{p-1}}\xi_2)-\frac13) (d\xi_2^2+dy_3^2), \\
H = -1, \;\; U d(\xi_2+i y_3)^3=-i\frac{\sqrt{3}(2p-\sqrt{3}i-1)^3}{576(p-1)^{3/2}} d(\xi_2+i y_3)^3.
\end{eqnarray*}
Substituting $g$, $H$, $U'=|U|$ into frame equations \eqref{eqframe2} and solving them, we can get the following theorem.
\bthm The associated family of the complete hyperbolic affine sphere in \cite{Roland Hildebrand.{2013}}:
\[
\mathbb{R}_{++} \times \mathbb{R}  \ni \bpm \xi \\  \mu \epm \rightarrow  \bpm r_1 \\  r_2 \\ r_3 \epm
=\left(1+\frac{\xi\sqrt{\xi+p+q-1}}{\sqrt{p+q}}\right)^{\frac{1}{3}} \bpm t^{-\frac{1}{3}}e^{\frac{q+1}{3q}\mu} \\ t^{-\frac{1}{3}}e^{-\frac{p+1}{3p}\mu} \\ t^{\frac{2}{3}}e^{-\frac{p-q}{3(p+q)}\mu} \epm,
\]
\begin{align*}
    t = & (\sqrt{\xi+p+q-1}+\sqrt{p+q}) \left( \frac{\xi}{(\sqrt{\xi+p+q-1}+\sqrt{p+q-1})^2}\right)^{\frac{\sqrt{p+q-1}}{\sqrt{p+q}}} \\
    &\cdot \left( \frac{\sqrt{\xi+p+q-1}+\sqrt{q-1}}{\xi+p} \right)^{\frac{1}{p}} \left( \frac{\sqrt{\xi+p+q-1}+\sqrt{p-1}}{\xi+q} \right)^{\frac{1}{q}}
\end{align*}
is given by

\begin{align*}
 &\vec{r}_{\l} = \sqrt{\frac{p^2-p+1}{3(p-1)}} \cdot \notag
 \\
 &\bpm
 (\coth(\sqrt{\frac{p^2-p+1}{4(p-1)}}\xi_2)-\frac{1}{\sqrt{3}}(\sqrt{3}\sin(t)+\cos(t))) e^{\sqrt{\frac{p^2-p+1}{12(p-1)}}[(\sqrt{3}\sin(t)+\cos(t))\xi_2+(\sqrt{3}\cos(t)-\sin(t))y_3]}\\
 (\coth(\sqrt{\frac{p^2-p+1}{4(p-1)}}\xi_2)-\frac{1}{\sqrt{3}}(-\sqrt{3}\sin(t)+\cos(t))) e^{\sqrt{\frac{p^2-p+1}{12(p-1)}}[(-\sqrt{3}\sin(t)+\cos(t))\xi_2+(-\sqrt{3}\cos(t)-\sin(t))y_3]}\\
 - \frac{\sqrt{p-1}}{p}(\coth(\sqrt{\frac{p^2-p+1}{4(p-1)}}\xi_2)+\frac{2}{\sqrt{3}}\cos(t)) e^{\sqrt{\frac{p^2-p+1}{3(p-1)}}(-\cos(t)\xi_2+\sin(t)y_3)}\\
 \epm,
 \end{align*}
 with the affine metric, the affine mean curvature and the affine cubic form:
 \begin{eqnarray*}
g = \frac {p^2-p+1}{4(p-1)}(\coth^2(\frac{\sqrt{p^2-p+1}}{2\sqrt{p-1}}\xi_2)-\frac13) (d\xi_2^2+dy_3^2), \\
H = -1, \;\; U d(\xi_2+i y_3)^3=e^{3it} \cdot \frac{\sqrt{3}}{72}\left(
\frac{p^2-p+1}{p-1}
\right)^{\frac{3}{2}}d(\xi_2+i y_3)^3.
\end{eqnarray*}
\ethm
$\mathbf{c=-2(p+q), \beta=+\infty.}$ In this case $s=-1$, and \eqref{eqaffsim1}-\eqref{eqaffsim3} the becomes
\begin{align*}
&\mathbb{R}_{--} \times \mathbb{R}  \ni \bpm \xi_2 \\  y_3 \epm \rightarrow  \vec{r} = \bpm r_1 \\  r_2 \\ r_3 \epm
 =\\
& \bpm
\frac{1}{\sqrt{3(p-1)}} (\sqrt{p^2-p+1} \coth(\frac{\sqrt{p^2-p+1}}{2\sqrt{p-1}}\xi_2)-1)e^{\frac{1}{2\sqrt{p-1}}\xi_2+\frac{2p-1}{2\sqrt{3(p-1)}}y_3} \\
\frac{1}{\sqrt{3(p-1)}} (\sqrt{p^2-p+1} \coth(\frac{\sqrt{p^2-p+1}}{2\sqrt{p-1}}\xi_2)-(p-1))e^{\frac{\sqrt{p-1}}{2}\xi_2-\frac{p+1}{2\sqrt{3(p-1)}}y_3} \\
-\frac{1}{\sqrt{3}p} (\sqrt{p^2-p+1} \coth(\frac{\sqrt{p^2-p+1}}{2\sqrt{p-1}}\xi_2)+p)e^{-\frac{p}{2\sqrt{p-1}}\xi_2-\frac{p-2}{2\sqrt{3(p-1)}}y_3}
\epm.
\end{align*}
This affine sphere is asymptotic to the cone given in case 3 of Theorem \ref{thm}. Since this affine sphere is in the same associated family with the previous case, their affine metrics and affine mean curvatures  are the same, and the cubic form of this case is
\[
U d(\xi_2+i y_3)^3=-i\frac{\sqrt{3}(2p+\sqrt{3}i-1)^3}{576(p-1)^{3/2}} d(\xi_2+i y_3)^3.
\]

\ms
\ms
At last, we explain that any regular convex cone in $\mathbb{R}^3$ should have a natural associated family of such cones.

If we try the transformations: $\sqrt{\frac{p^2-p+1}{12(p-1)}}\xi_2 \rightarrow \xi_2$, $\sqrt{\frac{p^2-p+1}{12(p-1)}}y_3\rightarrow y_3$, and
\[
\vec{r}_{\l} \rightarrow A \cdot (D \cdot B)^{-1} \cdot E \cdot \vec{r}_{\l}
\]
for $\vec{r}_{\l}$ in the Theorem 3.3, with the matrices A, B, D in the Theorem 3.1 and
\[
E = \left(
       \begin{array}{ccc}
         0 & 0 & \sqrt{\frac{3p^2}{p^2-p+1}} \ms \\
         0 & -\sqrt{\frac{3(p-1)}{p^2-p+1}} & 0 \ms \\
         -\sqrt{\frac{3(p-1)}{p^2-p+1}} & 0 & 0\\
       \end{array}
     \right),
\]
it is easy to see that the associated family in Theorem 3.3 is the same with it in Theorem 3.1, i.e., case 1, case 3 and case 5 are all in the same associated family. Conversely, any member in this family can be found in these three cases. The cones which these members are asymptotic to also can be found in the corresponding cases of the semi-homogeneous cones, which form a natural associated family of cones.

There are similar results for case 4. To show this, we should use the following relationship between Weierstrass  $\mathscr{P}$-function and Jacobi SN-function:
\begin{equation}
    \mathscr{P} (x,g_2,g_3) = \frac{e_1-e_3}{sn^2(\sqrt{e_1-e_3}x,k^2)}+e_3,
\end{equation}
here $e_i$ are the roots of $4z^3-g_2z-g_3=0$ which are ordered by $e_1>e_2>e_3$, and $sn$ is the Jacobi SN-function with the  elliptic modulus $k^2 = \frac{e_2-e_3}{e_1-e_3}$.
For any associated family in Theorem 3.2 with $p=p_1$ and $c=c_1$, we try the transformation:
\[
\sqrt{e_1-e_3} \xi_2 \rightarrow \xi_2, ~~~~~\sqrt{e_1-e_3} y_3 \rightarrow y_3.
\]
For all $p \in [2,+\infty]$, we can get the parameter $c$ by solving the equation:
\[
k^2(p_1,c_1) = k^2(p,c).
\]
For example, if we choose $p_1 = 3$, $c_1 = 1$ and $p = 2$, we can get $c = 4.39$. In this way, we can find out $g$ and the module of $U$ are the same for $(p_1,c_1)$ case and $(p,c)$ case, i.e., they are in the same associated family. So we find out the corresponding relationship between the complete hyperbolic affine spheres in case 4 and the associated families in Theorem 3.2. By this relationship, we can find the corresponding cones of the associated families in Theorem 3.2, which also form natural  associated families of cones.

In conclusion, in general any regular convex cone in $\R^3$ has a natural associated $S^1$-family of such cones, which deserve further studies. Given any natural associated family of some semi-homogeneous cone, it may contain case 1, case 3 and  case 5 of semi-homogenous cones or case 4.
\section*{Acknowledgments}

The second author would like to express the deepest gratitude for the support of the Hong Kong University of Science and Technology during the project, especially Min Yan, Yong-Chang Zhu, Bei-Fang Chen and Guo-Wu Meng. The authors had also been supported by the NSF of China (Grant Nos. 10941002, 11001262), and the Starting Fund for Distinguished Young Scholars of Wuhan Institute of Physics and Mathematics (Grant No. O9S6031001).


\begin{thebibliography}{99}

\bibitem{Eugenio Calabi.{1972}}
Calabi, E., \emph{Complete affine hyperspheres I}, In Symposia Mathematica Vol. 10, pages 19-38.
Istituto Nazionale di Alta Matematica, Acad. Press, 1972.
%

\bibitem{Cheng-Yau.1980}
Cheng, S.Y. and Yau, S.T., \emph{On the existence of a complete K\"ahler metric on noncompact
complex manifolds and the regularity of Fefferman¡¯s equation}, Comm. Pure Appl. Math.,
33(4):507-544, 1980.
%

\bibitem{Cheng-Yau.1982}
Cheng, S.Y. and Yau, S.T., \emph{The real Monge-Amp\`{e}re equation and affine flat structures},
In Proc. 1980 Beijing Symp. on Diff. Geom. and Diff. Eq., Vol. 1, pages 339-370, 1982.
%

\bibitem{Cheng-Yau.1986}
Cheng, S.-Y. and Yau, S.-T.,
\emph{{C}omplete affine hyperspheres. part I. The completeness of affine metrics},
Communications on Pure and Applied Mathematics, \textbf{39}(6):839-866, 1986.
%

\bibitem{Chern.1955}
Chern, S.S., \emph{An elementary proof of the existence of isothermal parameters on a surface}, Proc. Amer. Math. Soc. \textbf{6} (5): 771¨C782, doi:10.2307/2032933, JSTOR 2032933.
%

\bibitem{DoEi01}
Dorfmeister, J. and Eitner, U., \emph{{W}eierstrass-type representation of affine sphere},  Abh. Math. Sem. Univ. Hamburg, \textbf{71} (2001), 225--250.
%


\bibitem{Dor-Wang01}
Dorfmeister, J. and Wang, E.X., \emph{{D}efinite affine spheres via loop groups I: general theory},  preprint.
%

\bibitem{Dor-Wang02}
Dorfmeister, J. and Wang, E.X., \emph{{D}efinite affine spheres via loop groups II: equivariant solutions},  draft.
%

\bibitem{DuPl09}
Dunajski, M. and Plansangkate, P., \emph{{S}trominger-Yau-Zaslow geometry, affine spheres and Painlev\'e  III}, Commun. Math. Phys., \textbf{290} (2009), 997--1-24.
%

\bibitem{Fox.2012}
Fox, D.J.F., \emph{A Schwarz lemma for K\"ahler affine metrics and the canonical potential of a
proper convex cone}, arXiv e-print math.DG:1206.3176, 2012.
%

\bibitem{Gue79}
Guest, M., \emph{Harmonic Maps, Loop Groups, and Integrable Systems}, Cambridge University Press, 1997
%

\bibitem{Roland Hildebrand.2012}
Hildebrand, R., \emph{Einstein-Hessian barriers on convex cones}. Optimization Online e-print
2012/05/3474, 2012. Accepted at Math. Oper. Res.
%

\bibitem{Roland Hildebrand.{2013}}
Hildebrand, R., \emph{Analytic formulas for complete hyperbolic affine spheres}, Beitrage zur Algebra und Geometrie/Contributions to Algebra and Geometry, (2013), 1-24.
%


\bibitem{Kap97}
Kaptsov, O. V. and Shan'ko, Yu. V., \emph{Trilinear representation and the Moutard transformation for the Tzitzeica equation}, arXiv:solv-int/9704014v1
%

\bibitem{Li-Simon-Zhao.{1993}}
Li, A.M., Simon, U. and Zhao, G.S., \emph{Global affine differential geometry of hypersurfaces},
volume 11 of De Gruyter expositions in mathematics. Walter de Gruyter, 1993.
%

\bibitem{Lin-Wang-Wang}
Lin, Z.C., Wang, G. and Wang, E.X. \emph{Dressing actions on proper definite affine spheres}, arXiv:1502.04766.
%

\bibitem{Loftin.2002}
Loftin, J.C., \emph{Affine spheres and K\"ahler-Einstein metrics}, Math. Res. Lett., 9(4):425-432, 2002.
%

\bibitem{LoYauZa05}
Loftin, J., Yau, S.-T. and Zaslow, E.,
\emph{{A}ffine manifolds, SYZ geometry and the Y vertex},
J. Differential Geom., \textbf{71}(1):129-158, 2005.
%

\bibitem{Simon-Wang.1993}
Simon, U. and Wang, C.P., \emph{Local theory of affine 2-spheres}. Proc. Symposia Pure
Math. 54 (1993), 585-598.
%

\bibitem{Spivak}
Spivak, M., \emph{A Comprehensive Introduction to Differential Geometry 4 (3 ed.)}. Publish or Perish, 314-346


\bibitem{Te08}
Terng, C.L., \emph{{G}eometries and symmetries of soliton equations and integrable elliptic equations}, Surveys on geometry and integrable systems, 401--488, Adv. Stud. Pure. Math., 51, Math. Soc. Japan, Tokyo, 2008.
%

\bibitem{Uh89}
Uhlenbeck, K., \emph{{H}armonic maps into Lie groups (classical solutions of the chiral model)}, J. Diff. Geom., \textbf{30} (1989), 1--50.
%

\bibitem{Wa06}
Wang, E.X., \emph{{T}zitz\'{e}ica transformation is a dressing action}, J. Math. Phys., \textbf{47} (2006), no. 5,  053502, 13 pp.
%

\bibitem{Wang-lin-Wang}
Wang, G., Lin, Z.C. and Wang, E.X., \emph{Permutability theorem for definite affine spheres and the group structure of dressing actions.}, preprint.
%

\bibitem{Zak79}
Zakharov, V.E. and Shabat, A.B., \emph{Integration of non-linear equations of mathematical
physics by the inverse scattering method, II},Funct. Anal. Appl., \textbf{13} (1979), 166-174.
%

\end{thebibliography}
\end{document}